\documentclass{OAGM}
\OAGMarXiv{1404.3538}
\usepackage{graphics}
\usepackage[utf8]{inputenc}
\usepackage{lmodern}
\usepackage{hyperref}
\usepackage{amsmath}
\usepackage{amsthm}
\usepackage{amsfonts}
\usepackage{amssymb}
\usepackage{graphicx}
\usepackage{caption}
\usepackage[labelformat=simple]{subcaption}
\usepackage{algorithm}
\usepackage{algorithmic}

\theoremstyle{plain}

\theoremstyle{definition}

\providecommand{\abs}[1]{\lvert#1\rvert}
\providecommand{\scal}[2]{\langle #1,#2 \rangle}

\providecommand{\norm}[1]{\lVert #1 \rVert}
\providecommand{\ndel}[1]{\left(#1\right)}

\providecommand{\ie}{i.\,e.\ }
\providecommand{\eg}{e.\,g.\ }

\title{Sinogram constrained TV-minimization for metal artifact reduction in CT }

\author{Clemens Schiffer%
  \thanks{Institute for Mathematics and Scientific Computing, University of Graz, Heinrichstraße 36, A-8010 Graz, Austria
  \texttt{(clemens.schiffer@edu.uni-graz.at)}}
  and Kristian Bredies%
  \thanks{Institute for Mathematics and Scientific Computing, University of Graz, Heinrichstraße 36, A-8010 Graz, Austria
  \texttt{(kristian.bredies@uni-graz.at)}}}

\begin{document}
\maketitle

\begin{abstract}
  A new method for reducing metal artifacts in X-ray computed tomography (CT) images
  is presented. 
  It bases on the solution of a convex optimization problem with inequality constraints on the sinogram,
  and total variation regularization for the reconstructed image.
  The Chambolle-Pock algorithm is used to numerically solve the discretized version of the optimization problem.
  As proof of concept we present and discuss numerical results for synthetic data.
\end{abstract}
%
%
\section{Introduction} \label{sec:Intro}
X-ray computed tomography (CT) has become a major tool in medical imaging, 
yet metal inclusions in the patient, such as pacemakers, screws, orthopedic implants and dental fillings, 
still pose a problem. 
The basic principle of CT is rotating an X-ray source and an opposing detector 
around the object to be scanned, producing projection data from a range of angles.
These projections of intensities, known as the sinogram, have to be reconstructed into a spatial image,
in order to obtain the desired cross-section.
To achieve this, the standard method is filtered back projection (FBP), which works very well under ideal or even normal circumstances. 
Still errors from various sources, such as patient movement, beam scatter etc.~disturb this reconstruction process 
and lead to artifacts. 
In particular, the presence of metal inclusions produces significant dark and bright smears,
known as metal artifacts, that may render the image unusable for diagnostics.


Existing techniques for reducing metal artifacts can be categorized into technical enhancements and computational methods.
The former include the use of dual energy CT devices~\cite{Bamberg_dual}, 
which feature two sources that produce projections from different angles. The different energy levels 
used by the two sources allow to identify a substance by their specific absorption spectrum,
this information can then be used to reduce metal artifacts.
In contrast, computational methods focus on processing the already collected data in a way that removes the artifacts.
This is, for instance, accomplished by replacing the corrupted data either in the spatial or the projection domain,
e.\,g., by interpolation~\cite{Zhang_interp_fracO} or inpainting~\cite{Chen_inpainting},
with various refinements such as repeated projection, back-projection and replacement of bad data \cite{Boas_2iter,karimi2012segmentation}. 
Variational methods involving total variation regularization have been proposed~\cite{Ritschl_tv}, 
and applied with the recently developed Chambolle-Pock algorithm ~\cite{sid_tv_cp}.
The novelty of our method lies in the dedicated treatment of metal artifacts by introducing pointwise 
inequality constraints on the sinogram in addition to total variation regularization of the image.


In Section \ref{sec:ProbMod} the problem is presented and mathematically modeled,
in Section \ref{sec_alg} the numerical solution is discussed and an algorithm for solving the proposed optimization problem is described, in Section \ref{sec_exp}
results are presented and, finally, in Section \ref{sec_concl}, a conclusion is given.    
      
%
%
\section{Problem and Mathematical Model} \label{sec:ProbMod}
The object to be scanned is modeled as a continuous function $f(x,y)$ in $\mathbb{R}^2$, that represents the density  of the 
object at a point $(x,y)$. It is assumed to be zero outside a region of interest.
The assumption of the model is that the attenuation an X-ray suffers, \ie the loss of intensity, is proportional to the traversed density.  
Therefore, the process of an X-ray passing through the object, being weakened and received by a detector is modeled as 
integrating $f$ over a straight line $L$ 

\begin{align}
(Af)(\varphi, s) = \int_L f(x) ds,
\end{align} 
where $\varphi$ is the gradient angle of $L$ and $s$ its signed distance from the origin, the \emph{offset}.
This is called the \emph{Radon transform} and the data it produces a \emph{sinogram}, which
is commonly interpreted as a function of angle and offset.  
Consider a ray that starts at the source with intensity $I(x_S)$ and whose intensity is measured at the detector giving $I(x_D)$.
By the modeling assumption~\cite{Kak_princ_ct} we get

\begin{align}
(Af)(\varphi, s) = \int_L f(x) ds = -\log(I(x_D)) + \log(I(x_S)).
\end{align}
Now if the ray passes through a very dense area, and is attenuated so strongly that the detector cannot 
differentiate its signal from noise,
\ie the received signal is unknown but lower than a certain threshold $\varepsilon$, then $I(x_D) < \varepsilon$.
The intensity at the source  $I(x_S)$ is considered constant, and therefore, we get
$(Af)(\varphi, s) > -\log(\varepsilon) + \log(I(x_S)) = \colon C$.
Thus, this process of beam cancellation is represented by ``capping'' the sinogram, 
\ie defining a certain threshold $C$ and considering the values for each point in the sinogram that exceeds $C$ 
as unknown, but at least as great as the threshold.
The domain of the sinogram is denoted as $\Omega$, while the portion where the threshold is exceeded is denoted as $\Omega_0$.
The capped sinogram is called $U_0$ and considered to be the given data.
It is assumed - as in total variation denoising - that $f$ has a certain spatial structure for the propose of this paper it is assumed
to be piecewise constant and therefore to admit a low total variation semi-norm. 
In trying to reconstruct $f$ we therefore look for an approximation $u$ that has also low TV-norm and
for which an application of the Radon transform will produce a sinogram that is similar to that of $f$ in $\Omega \setminus \Omega_0$
and also has values greater or equal than $C$ in $\Omega_0$. 
This leads to the optimization problem

\begin{align} \label{opt_main}
\begin{cases}
\min_{ u \in X  } \quad \frac{1}{2} \norm{ Au - U_0}^2 
+ \lambda \norm{ \nabla u}_1 \\
\text{s.t. }  Au\mid_{\Omega_0} \geq C
\end{cases}
\end{align} 
where
$\norm{ \cdot} =\norm{ \cdot}_{L^2(\Omega \setminus \Omega_0)}$ and $ 
\norm{ \cdot}_1 =\norm{ \cdot}^2_{L^1(\Omega) }$.
Here $\lambda > 0$ serves as a balancing parameter between regularization via the total variation semi-norm, and the discrepancy 
between $Au$ and $U_0$, which is measured on $\Omega \setminus \Omega_0$ \ie only on sinogram data that is considered to be correct.   
Using the indicator function

\begin{align}
I_M(x) 
=  \begin{cases}
	0,      &\text{if } x \in M  \\
	\infty, &\text{if } x \notin M \\
	\end{cases}
\end{align}
where $0 \cdot \infty = 0$ and $ 1 \cdot \infty = \infty$ are set,
we can reformulate \eqref{opt_main} as

\begin{align}	\label{opt_main_unconstr}
	\min_{u\in X} \quad \frac{1}{2} \norm{ Au - U_0}^2 + \lambda \norm{ \nabla u}_1 + I_{ \{Au\mid_{\Omega_0} \geq C \}}(u). 
\end{align} 
%
%
%
The parameter $\lambda$ is not needed if one enforces $ Au = U_0$ in $\Omega \setminus \Omega_0$ as a hard constraint.
Then, the minimization problem becomes

\begin{align} \label{opt_main_hc}
	\begin{cases}
	\min_{ \substack {u\in X } } \quad  \norm{ \nabla u}_1 \\
	\text{s.t. } Au\mid_{\Omega_0} \geq C \wedge Au\mid_{\Omega \setminus \Omega_0} = U_0
	\end{cases}
\end{align}
which is equivalent to

\begin{align} \label{opt_main_hc_unconstr}
	\min_{ \substack {u\in X } } \quad  \norm{ \nabla u}_1 
	+  I_{ \{Au\mid_{\Omega_0} \geq C \}}(u) + I_{ \{Au\mid_{\Omega \setminus \Omega_0} = U_0\}}(u).
\end{align}

%
%
\begin{figure}
        \centering
        \begin{subfigure}[b]{0.3\textwidth}
                \includegraphics[width=\textwidth]{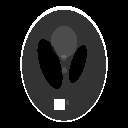}
                \caption{}
                \label{fig:gt}
        \end{subfigure}   \raisebox{5.5\height}{\scalebox{3}{$\rightarrow$}}
        \begin{subfigure}[b]{0.5\textwidth}
                \includegraphics[width=\textwidth]{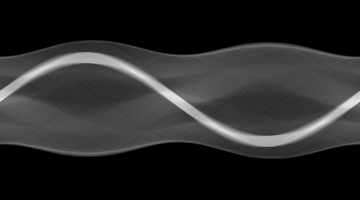}
                \caption{}
                \label{fig:sino_uncap}
        \end{subfigure}%
        
        \begin{subfigure}[b]{0.3\textwidth}
                \includegraphics[width=\textwidth]{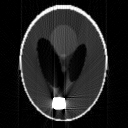}
                \caption{}
                \label{fig:bad}
        \end{subfigure}
        \raisebox{5.5\height}{\scalebox{3}{$\leftarrow$}}
        \begin{subfigure}[b]{0.5\textwidth}
                \includegraphics[width=\textwidth]{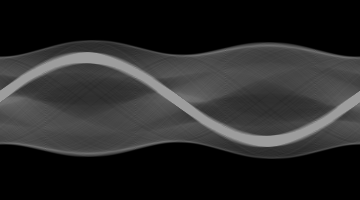}
                \caption{}
                \label{fig:sino_capped}
        \end{subfigure}
	    \caption{Metal artifacts are simulated with a synthetic test image. \subref{fig:gt} shows a Shepp-Logan phantom,
		with values in $[0,1]$ to which a block of value 3 has been added, simulating metal.
		\subref{fig:sino_uncap} shows the result of applying the Radon transform, \ie the corresponding sinogram. 
		In \subref{fig:sino_capped}, this sinogram has been capped: values greater than $45$ have been set to $45$,
		\subref{fig:bad} shows the devastating results of simply applying the MATLAB-function \emph{iradon} to the capped sinogram.
		Figures \subref{fig:gt} and \subref{fig:bad} display only values in  $[0,1]$ for better contrast.
	      }\label{fig:synPics} 
\end{figure}

%
%
\section{Numerical Solution} \label{sec_alg}
We are numerically solving a discretized version of \eqref{opt_main_unconstr} and \eqref{opt_main_hc_unconstr}, following the presentation in \cite{KB}.
In order to discretize the gradient $\nabla$ a finite difference scheme with spacing $h$ and zero boundary extension is used. Note that then, 
$\norm{\nabla_h}^2 < \frac{8}{h^2}$ and its adjoint is $ \nabla_h^\ast = - \operatorname{div}_h$ with a similar discretization.
The discretization of the Radon transform $A_h$ is implemented as follows:
Every point in the image domain is, for each of the $N$ angles, projected onto the detector line, which is 
separated into $M$ bins. The value at the point is then assigned proportionally to the two bins the point is projected inbetween.
This implementation yields the same results as the MATLAB-function \emph{radon}, when using $2\times2$ oversampling.
The adjoint of the discrete Radon transform $A_h^\ast$ is a discretization of the linear back projection, which is implemented similarly,
each point in the image domain is projected onto the detector line and linear interpolation is performed. Then the values for all angles are summed up.
Numerical tests confirm that this implementation provides the adjoint to the discrete Radon transform. 
 
The optimization problems \eqref{opt_main_unconstr} and \eqref{opt_main_hc_unconstr} are of the form 

\begin{equation} \label{opt_abst}
\min_{u\in X} F(u) + G(Ku),
\end{equation}
with $F \colon X \rightarrow \mathbb{R}, G\colon Y \rightarrow \mathbb{R}$ convex, lower semi-continuous and proper, as well
as $K \colon X \rightarrow Y$ linear and continuous. 
These problems satisfy the sufficient conditions for the Fenchel-Rockafellar duality~\cite{KB}.
The dual problem reads as 

\begin{align}
\max_{w \in Y^\ast}-F^\ast(-K^\ast w) - G^\ast(w). 
\end{align} 
Now, solving the primal and dual equation simultaneously can be interpreted as 
finding the saddle-point of the function

\begin{align} \label{L_saddle}
L(u,w) = \scal{w}{Ku} + F(u) - G^\ast(w).
\end{align}
In order to solve this problem, the Chambolle-Pock algorithm \cite{Pock_dual}

\begin{align} \label{AH-Iter}
	\begin{cases}
	w^{k+1} = (id + \tau   \partial G^\ast)^{-1}(w^{k} + \tau   K \bar{u}^k )   \\
	u^{k+1} = (id + \sigma \partial F     )^{-1}(u^{k} - \sigma K^\ast w^{k+1}) \\ 
	\bar{u}^{n+1} = 2u^{k+1} - u^k
	\end{cases}
\end{align}
is employed, which converges to a saddle point of \eqref{L_saddle}, if $\sigma \tau \norm{K}^2 < 1$.
In order to match the form of \eqref{opt_abst}, the following discrete version of \eqref{opt_main}

\begin{align} \label{opt_main_h}
\min_{ u \in \mathbb{R}^{n\times m}} \quad \frac{1}{2} \norm{ A_h u - U_0}^2 + \lambda \norm{ \nabla_h u}_1 + I_{ \{A_h u\mid_{\Omega_0} \geq C\}}(u)
\end{align}
is dualized with $F = 0$ and $G \colon \mathbb{R}^{N \times M} \times \mathbb{R}^{n \times m \times 2} \rightarrow \mathbb{R}_\infty$ 
as well as  the linear mapping  $K \colon \mathbb{R}^{m \times n }  \rightarrow \mathbb{R}^{M \times N} \times \mathbb{R}^{m \times n \times 2}$ defined by

\begin{align}
	G(x,y) = \frac{1}{2}\norm{x - U_0}^2 +  I_{ \{x\mid_{\Omega_0} \geq C \}}(x) + \lambda \norm{y}_1 \qquad K = 
	\begin{bmatrix}
	A_h \\ \nabla_h
	\end{bmatrix}.
\end{align}
For $F = 0$ the resolvent is $ (id + \sigma \partial F)^{-1} = id$, as for $G^\ast$ its resolvent can be evaluated componentwise, \ie 
for $G(x,y) = G_1(x) + G_2(y)$ we have 

\begin{align}
(id + \sigma \partial G^\ast)^{-1}(\bar v, \bar w) = 
    \begin{pmatrix}
	(id + \sigma \partial G^\ast_1)^{-1}(\bar v) \\
	(id + \sigma \partial G^\ast_2)^{-1}(\bar w)
	\end{pmatrix}
\end{align}
As $G_1(x) = \sum_{(i,j)\in \Omega \setminus \Omega_0} (x-U_0)_{i,j}^2 + \sum_{(i,j) \in \Omega_0} I_{x_ \geq C}(x_{i,j}) $, its dual is

\begin{align}
G_1^\ast(\xi) = \sum_{(i,j) \in \Omega \setminus \Omega_0 } g_1^\ast (\xi_{i,j}) +  \sum_{(i,j) \in  \Omega_0 } \bar g_1^\ast (\xi_{i,j})
\end{align}
with $g_1^\ast (\xi_{i,j}) = 	\frac{1}{2}\xi_{i,j}^2 + \xi_{i,j}(U_0)_{i,j}$
and
 
\begin{align}
 \bar g_1^\ast (\xi_{i,j}) = 
 \begin{cases}	\infty, &\text{if } \xi_{i,j} > 0  \\
	\xi_{i,j}C,                 &\text{if } \xi_{i,j} \leq 0 \\
 \end{cases}
\end{align}
consequently the resolvent can be calculated to be

\begin{align}
(id + \tau \partial G^\ast_1)^{-1}(\bar v_{i,j}) = 
	\begin{cases}
	\frac{\bar v_{i,j}-\tau (U_0)_{i,j}}{1+\tau},     &\text{if } (i,j) \in \Omega \setminus \Omega_0   \\
	\min \{ \bar v_{i,j} - \sigma C, 0 \},                 &\text{if } (i,j) \in  \Omega_0 .
	\end{cases}
\end{align}
The dual of $G_2(y) = \lambda \norm{y}_1$ is $G_2^\ast(\eta) = I_{B_\lambda^\infty(0)}(\eta)$, where
${B_\lambda^\infty(0)}$  denotes the norm ball of radius $\lambda$ around the origin 
in the maximum norm, consequently the resolvent of $G_2^\ast$ is the projection onto ${B_\lambda^\infty(0)}$:

\begin{align}
\ndel{(id + \tau \partial G^\ast_2)^{-1}(\bar w)}_{i,j}
 = \ndel{ \mathcal{P}_{B_\lambda^\infty(0)}(\bar w) }_{i,j} 
 = \frac{\bar w_{i,j}}{ \max \{ 1, \abs{ \bar w_{i,j} } / \lambda \} }.
\end{align} 
Then, the Chambolle-Pock iteration provides us with Algorithm \ref{AH-alg}.

\begin{algorithm} 
\caption{Chambolle-Pock algorithm for minimizing \eqref{opt_main_h}}
\begin{algorithmic} \label{AH-alg}
\REQUIRE  \textbf{Input}:  $U_0 \in \mathbb{R}^{N \times M}, C >0, \, \sigma >0, \, \tau >0, \, h>0 $
\REQUIRE $\sigma \tau < (\norm{A}^2 + \frac{8}{h^2})^{-1}$
\STATE $\Omega_0 \leftarrow \{ (i,j) \mid (U_0)_{i,j} \geq C \}$
\FOR{$k < k_{max}$} 
\STATE $\bar w \leftarrow w + \tau\nabla_h \bar u$
\STATE $ w \leftarrow \mathcal{P}_{B_\lambda^\infty(0)}(\bar w)$
\STATE $\bar v \leftarrow v + \tau A \bar u$
\STATE $v_{i,j} \leftarrow 
\begin{cases}
\min \{ \bar v_{i,j} - \sigma C, 0 \}, &\text{if } (i,j) \in \Omega_0  \\
\ndel{\bar v_{i,j} - \tau (U_0)_{i,j}}/(1+\tau), &\text{else } 
\end{cases}
$
\STATE $ u_{new} \leftarrow u + \sigma div_h(w) -\sigma A^\ast(v) $
\STATE $ \bar u \leftarrow 2 u_{new} - u$
\STATE $u \leftarrow u_{new}$
\ENDFOR
\end{algorithmic}
\end{algorithm}
In the case of hard constraints we have 
$G_1 =  I_{\{x\mid_{\Omega_0} \geq C\}} + I_{\{x\mid_{\Omega \setminus \Omega_0} = U_0\}} $
and the resolvent of $G^\ast_1$ is again
$(id + \sigma \partial G_1^\ast)^{-1}(\bar v_{i,j}) = \min \{ \bar v_{i,j} - \sigma C, 0 \}$ in $\Omega_0$ and
in $\Omega \setminus \Omega_0$ we have  $(id + \sigma \partial G_1^\ast)^{-1}(\bar v_{i,j}) = \bar v_{i,j} - \sigma (U_0)_{i,j}$
finally we can choose $\lambda>0$ arbitrarily  \eg $\lambda = 1$.
%
%
\begin{figure}
        \centering
        \begin{subfigure}[b]{0.4\textwidth}
                \includegraphics[width=\textwidth]{syn_gt.png}
                \caption{Ground truth}
                \label{fig:gt2}
        \end{subfigure}
        \begin{subfigure}[b]{0.4\textwidth}
                \includegraphics[width=\textwidth]{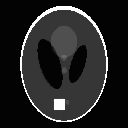}
                \caption{Result hard constraint}
                \label{fig:result_hc}
        \end{subfigure}
        
        \begin{subfigure}[b]{0.3\textwidth}
                \includegraphics[width=\textwidth]{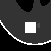}
                \caption{}
                \label{fig:gt_zoom}
        \end{subfigure}
        \begin{subfigure}[b]{0.3\textwidth}
                \includegraphics[width=\textwidth]{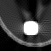}
                \caption{}
                \label{fig:bad_zoom}
        \end{subfigure} 
               \begin{subfigure}[b]{0.3\textwidth}
                \includegraphics[width=\textwidth]{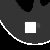}
                \caption{}
                \label{fig:res_hc_zoom}
        \end{subfigure} 
         \label{fig:synPics2}
        \caption{Results: \subref{fig:gt2} shows the ground truth again, \subref{fig:result_hc} shows the result of the proposed method 
        with hard constraints after 80000 iterations with Peak signal-to-noise ratio PSNR $ = 47.6$\,dB.
        In the second row details are shown, note the little area of tissue next to the metal in the ground truth \subref{fig:gt_zoom},
        which along with other details is completely obscured in the reconstruction of the capped sinogram \subref{fig:bad_zoom}.
        In our result with hard constraints \subref{fig:res_hc_zoom} the tissue next to the metal can at least be recognized,
        other details are preserved, but note the general blurriness.    } 
\end{figure}
%
%
\begin{figure}
        \centering
        \begin{subfigure}[b]{0.4\textwidth}
                \includegraphics[width=\textwidth]{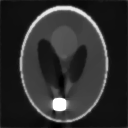}
                \caption{Result without constraints}
                \label{fig:res_unconst}
        \end{subfigure}
        \begin{subfigure}[b]{0.4\textwidth}
                \includegraphics[width=\textwidth]{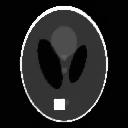}
                \caption{Result with constraints}
                \label{fig:result_l2}
        \end{subfigure}
        \caption{Results for reconstruction from a capped sinogram with 5\% Gaussian noise,
        \subref{fig:res_unconst} shows the result of L2-TV minimization without constraints similar to the method presented
        in~\cite{sid_tv_cp} if one does not account for the presence of metal, \subref{fig:result_l2} shows the
        result of the proposed method of L2-TV minimization with constraints after 80000 iterations with Peak
        signal-to-noise ratio PSNR $ = 40.1$\,dB. The parameter $\lambda = 10^{-4.1}$ has been	 
		manually chosen, to optimize PSNR.} 
\end{figure}
%
%
\section{Experiments and Results} \label{sec_exp}
In order to test the proposed method, a Shepp-Logan phantom of size $128 \times 128$ pixels and values in $[0,1]$ was 
used. In an area of $10 \times 10$ pixels the value $3$ was added to representing a metallic object, see Figure \ref{fig:gt}.
Applying the Radon transform results in the sinogram in Figure \ref{fig:sino_uncap}.
Simply trying to reconstruct from the capped sinogram \eg using the MATLAB\textsuperscript{\textregistered} function \emph{iradon} will result in 
artifacts (see Figure \ref{fig:bad})  that are very similar to those seen in real reconstructions.
Algorithm \ref{AH-alg} was applied to the capped sinogram $U_0$ with the operator $A = \frac{1}{D}A_h$ where $D$ is a crude bound 
for the norm of $A_h$. Then $\norm{A} < 1$ and consequently, we used $\sigma = \tau  < (1 + \frac{8}{h^2})^{-\frac{1}{2}}$, which led to a convergent algorithm.
The result for hard constraints after $80000$ iterations can be seen in Figure \ref{fig:result_hc}. 
Results for reconstruction the a capped sinogram with 5\%  Gaussian noise, using a $L^2$ data term after $80000$ iterations can be seen in Figure \ref{fig:result_l2}, the parameter $\lambda$ was manually chosen to optimize PSNR.
\section{Discussion/Conclusion} \label{sec_concl}
A new method of reducing metal artifacts has been presented, and the results serve as a proof of concept,
in particular the hard constrained optimization with synthetic data succeeds in removing the artifacts while preserving
most of the details.  
The method is not limited to errors produced by metal, but could be extended to deal with any kind of corrupted or missing data.
The most obvious drawback of the proposed method is the need for many iterations to produce satisfying results, which makes the method currently very slow. In our implementation this is mitigated to some extent by GPU-based parallelization, one iteration taking about 45\,ms still leads to 80000 iteration on the toy problem taking about an hour, which prevents practical applications.
We think that additional preconditioning may be necessary in order to provide adequate speed. 
We are optimistic that the method will also perform on more detailed and natural images without producing additional artifacts,
and thus effectively and robustly remove metal artifacts in CT.  A detailed study with respect to real data will be subject of future research. 
%
%
\bibliography{refs}
\end{document}